\documentclass[dvipdfmx]{amsart}

\makeatletter
\@addtoreset{equation}{section}

\makeatother

\usepackage{amsfonts, amsmath, amssymb, amsthm}
\usepackage{url}
\usepackage[dvipdfmx]{graphicx}
\usepackage[dvipdfmx]{color}
\usepackage{amscd}
\usepackage[right]{lineno}
\usepackage{enumerate}

\newtheorem{theorem}{Theorem}[section]
\newtheorem{lemma}{Lemma}[section]
\newtheorem{proposition}{Proposition}[section]
\newtheorem{corollary}{Corollary}[section]
\theoremstyle{remark}

\newcommand{\teich}{\mathcal{T}}

\newcommand{\ext}{{\rm Ext}}

\newcommand{\charfunc}{\boldsymbol{1}}

\newcommand{\Bers}[1]{\mathcal{T}^B_{#1}}

\newcommand{\ThursM}{\mu_{Th}}

\newcommand{\pushThursMBers}{{\boldsymbol \mu}^B}

\newcommand{\hol}{H}

\newcommand{\ray}{{\boldsymbol r}}

\makeindex

\begin{document}

\title[Bounded pluriharmonic functions]{Bounded pluriharmonic functions and holomorphic functions on Teichm\"uller space}
\author{Hideki Miyachi}

\date{\today}
\address{School of Mathematics and Physics,
College of Science and Engineering,
Kanazawa University,
Kakuma-machi, Kanazawa,
Ishikawa, 920-1192, Japan
}
\email{miyachi@se.kanazawa-u.ac.jp}
\thanks{This work is partially supported by JSPS KAKENHI Grant Numbers
20H01800,
20K20519,
22H01125}

\begin{abstract}
In this paper, we discuss the boundary behavior of bounded pluriharmonic functions on the Teichm\"uller space. We will show a version of the Fatou theorem that every bounded pluriharmonic function admits the radial limits along the Teichm\"uller geodesic rays, and a version of the F. and M. Riesz theorem that the radial limit of a non-constant bounded holomorphic function is not constant on any non-null measurable set on the Bers boundary in terms of the pluriharmonic measure. As a corollary, we obtain the non-ergodicity of the action of the Torelli group for a closed surface of genus $g\ge 2$ on the space of projective measured foliations.
\end{abstract}

\maketitle

\section{Introduction}

\subsection{Background}
\label{subsec:background}
Let $\teich_{g,m}$ be the Teichm\"uller space of Riemann surfaces of analytically finite type $(g,m)$ with $2g-2+m>0$. 
The Teichm\"uller space $\teich_{g,m}$ admits a natural complex structure and a natural complete distance, called the \emph{Teichm\"uller distance}, inherited from quasiconformal deformations of Riemann surfaces. Under the complex structure, the Teichm\"uller distance coincides with the Kobayashi distance (cf. \cite{Royden_2:1971}). If we fix a base point $x_0=(M_0,f_0)\in \teich_{g,m}$, the Teichm\"uller space $\teich_{g,m}$ is embedded as a bounded domain in the space of bounded holomorphic quadratic differentials on the mirror to $M_0$. 
The image $\Bers{x_0}$ of the image and the boundary $\partial \Bers{x_0}$ are called the \emph{Bers slice} and  the \emph{Bers boundary} with basepoint $x_0$, respectively (cf. \S\ref{subsec:Bers-slice}).
The Bers boundary is originated from the study by L. Bers in \cite{Bers:1970}, and since then, it is studied by many mathematician. 
 It is conjectured that the Bers boundary is fractal and self-similar at the fixed point with respect to the action of the pseudo-Anosov mapping class (See \cite{Canary:2010} , \cite{Komori_Sugawa_Wada_Yamashita:2006} and \cite[Problem 7 in Chapter 10]{McMullen:1996}). 

\subsection{Purpose of the research}
This paper is a part of the study of holomorphic and harmonic functions on the Teichm\"uller space in terms of the Thurston theory to approach above mentioned conjectures.


J-P. Demailly \cite{Demailly:1987}  establishes the Poisson kernels and the pluriharmonic measures for bounded hyperconvex domains in the complex Euclidean space. 
S. Krushkal \cite{Krushkal:1991} shows that the Bers slice is hyperconvex. 
H. Shiga \cite{Shiga:1984} shows that the Bers slice is polynomially convex, and hence every holomorphic function on $\teich_{g,m}$ can be approximated by holomorphic functions with the Poisson integral presentations.

In \cite{Miyachi:2023}, the author proves the Poisson integral formula for pluriharmonic functions on $\teich_{g,m}$ which are continuous on the Bers compactification, and gives a characterization of the Poisson kernels and the pluriharmonic measures in the sense of Demailly for the Bers slice $\Bers{x_0}$. Actually, it is shown that the pluriharmonic measure coincides with the pushforward measure of the (normalized) Thurston measure on $\mathcal{PMF}$ via the natural parametrization of b-groups without $APT$ in terms of the ending laminations on the Bers boundary (cf. \eqref{eq:curve-complex-mininal-bers} and \S\ref{subsec:Thurston_measure}). We also observe in \cite{Miyachi:2023} a version of Schwarz's theorem in \cite{Schwarz:1872} which studies the behavior of the Poisson integral of integrable functions at boundary points where given integrable functions are continuous.

Table \ref{table:dictionary} is a dictionary which clarifies the meaning of the abstract objects in the function theory in terms of the moduli of Riemann surfaces. To apply cultivated researches in the function theory to the Teichm\"uller theory, it is needed to increase the entries of the dictionary. 
Our researches are also developed with applying essentially the sophisticated researches of the hyperbolic geometry and the Kleinian groups (cf. \S\ref{subsec:BoundarygroupswithoutAPT}).
For instance, from the dictionary, the set of projective classes of uniquely ergodic measured foliations has two faces.
Indeed, the set is recognized as a set of almost all directions in the infinitesimal sphere as well as a full measure set as a subset of the Bers boundary.
Their recognitions are connected via the ending lamination theorem by J. Brock, D. Canary and Y. Minsky \cite{Brock_Canary_Minsky:2012} (cf. Proposition \ref{prop:Teichmuller-limit}).

\begin{table}
\begin{tabular}{|c|c|}\hline
Upper half-plane $\mathbb{H}$ & Teichm\"uller space $\teich_{g,m}$ \\ \hline\hline
Harmonic function & Pluriharmonic function \\ \hline
Compactification $\overline{\mathbb{H}}$ &
\begin{tabular}{l} 
Bers compactification (CA \cite{Bers:1970}) \\
Gardiner-Masur compactification (TOP-EL \cite{Gardiner_Masur:1991}) \\
Thurston compactification (TOP \cite{Douady_Fathi_Fried_Laudenbach_Poenaru_Shub:1979})
\end{tabular}
\\ \hline
Ideal boundary $\partial \mathbb{H}$ & 
\begin{tabular}{l} 
Bers boundary (CA \cite{Bers:1970}) \\
Gardiner-Masur boundary (TOP-EL \cite{Gardiner_Masur:1991}) \\
Thurston boundary (TOP \cite{Douady_Fathi_Fried_Laudenbach_Poenaru_Shub:1979})
\end{tabular}
\\ \hline
Hyperbolic metric  & 
\begin{tabular}{l}
Kobayashi-Royden Finsler metric (CA \cite{Royden_1:1971}) \\
Teichm\"uller metric (EL \cite{Royden_2:1971}) 
\end{tabular}
\\ \hline
Hyperbolic distance  & 
\begin{tabular}{l}
Kobayashi distance (CA \cite{Kobayashi:1967}) \\
Teichm\"uller distance (EL \cite{Royden_2:1971})
\end{tabular}
\\ \hline
Green function & 
\begin{tabular}{l}
Pluricomplex Green function (CA \cite{Demailly:1987}, \cite{Klimek:1985}) \\
$\log \tanh$ of the Teichm\"uller distance (EL \cite{Krushkal:1992}, \cite{Miyachi:2019})
\end{tabular}
\\ \hline
\begin{tabular}{l}
Horofunctions \\
(Busemann functions)
\end{tabular} & log of extremal lengths (EL \cite{Liu_Su:2014}, \cite{Miyachi:2014}) \\ \hline
Poisson kernel & 
\begin{tabular}{l}
Poisson kernel (CA \cite{Demailly:1987}) \\
Ratio of extremal lengths (EL \cite{Miyachi:2023})
\end{tabular}
\\ \hline
Harmonic measure on $\partial \mathbb{H}$ & 
\begin{tabular}{l}
Pluriharmonic measure (CA \cite{Demailly:1987}) \\
Normalized Thurston measure on $\mathcal{PMF}$ (TOP \cite{Miyachi:2023})
\end{tabular}
\\ \hline
\end{tabular}
\caption{A dictionary : TOP, EL, and CA stand for Topological, Extremal Length geometrical, and Complex Analytical aspects in the Teichm\"uller theory. Extremal length functions are plurisubharmonic (cf. \cite{Liu_Su:2017} and \cite{Miyachi:2019}). The Gardiner-Masur compactification and boundary work as mediators between TOP and CA via EL (e.g. \cite{Miyachi:2008}, \cite{Miyachi:2013}, \cite{Miyachi:2014}, and \cite{Miyachi:2015}). The Teichm\"uller distance and the extremal lengths are also treated from the topological and combinatorial viewpoints with the geometry of the curve complex. See \cite{Masur_Minsky:1999}, \cite{Masur_Minsky:2000}, \cite{Rafi_1:2007}, \cite{Rafi_2:2007} and \cite{Lenzhen_Rafi:2011} for instance.}
\label{table:dictionary}
\end{table}

\subsection{Results}
Given the natural development of the function theory, one of our next tasks is to understand the boundary behavior of (pluri)harmonic or holomorphic functions on $\teich_{g,m}$.

Our main results deal with the radial limits for bounded pluriharmonic functions on the Teichm\"uller space $\teich_{g,m}$.
As the preceding results, P. Fatou \cite{Fatou:1906} observes the existence of non-tangential limit for bounded harmonic functions on the unit disk $\mathbb{D}$ in $\mathbb{C}$ (cf. \cite[VII, \S3]{Nevanlina:1970} and \cite[Theorem IV.6]{Tsuji:1959}. See also \cite[Theorem 3.8.11]{Helms:2014} for the general case).
A. Kor{\'a}nyi \cite{Koranyi:1969} observes that bounded harmonic functions (in terms of the Bergman metric) admits the admissible limits (in the Kor{\'a}nyi sense) for the unit ball. E.M. Stein \cite{Stein:1972} discusses for the strictly pseudoconvex domains.
The radial limits are mostly considered for holomorphic (pluriharmonic) functions on geometrically nice domains, for instance for convex domains (e.g. \cite{Bagemihl_Seidel:1954}, \cite{Hakim_Sibony:1987}). For general domains, the formulation of ``radial" seems to be a delicate issue.

We say that a function $u$ on $\teich_{g,m}$ \emph{has the radial limit} if there are a measurable function $u^*$ on $\partial\Bers{x_0}$ with respect to the pluriharmonic measure and a full-measure set $\mathcal{E}_0$ on the space $\mathcal{PMF}$ of projective measured foliations with respect to the Thurston measure such that 
for any $x\in \teich_{g,m}$ and for any $[F]\in \mathcal{E}_0$, the limit of $u$ along the Teichm\"uller ray associated to the Hubbard-Masur differential for $F$ on $x$ exists, and coincides with $u^*$ at the limit point of the Teichm\"uller ray in $\partial \Bers{x_0}$.
We call the measurable function $u^*$ the \emph{radial limit} of $u$. Notice that the radial limit $u^*$ is assumed to be independent of the choice of the base point $x$.


\begin{theorem}[Radial limit]
\label{thm:main_radial_limit}
Any bounded pluriharmonic function $u$ has the radial limit almost everywhere on $\partial \Bers{x_0}$ with respect to the pluriharmonic measure, and the radial limit function $u^*$ is in $L^\infty(\partial \Bers{x_0})$.
\end{theorem}

The precise statement of Theorem \ref{thm:main_radial_limit} can be found in Theorem \ref{thm:radial_limit_theorem} in \S\ref{sec:radial_limit_theorem}.
From Theorem \ref{thm:main_radial_limit}, any bounded holomorphic function $f$ on $\teich_{g,m}$ has the radial limit $f^*\in L^\infty(\partial \Bers{x_0})$.
We also show a version of the F. and M. Riesz theorem for bounded holomorphic functions on the Teichm\"uller space as follows.

\begin{theorem}[Identity theorem]
\label{thm:F-M-Riesz}
A bounded holomorphic function $f$ on $\teich_{g,m}$ is constant if the radial limit $f^*$ of $f$ is constant on a non-null measurable set in $\partial\Bers{x_0}$ with respect to the pluriharmonic measure.
\end{theorem}
As a corollary to Theorem \ref{thm:F-M-Riesz}, we obtain
\begin{corollary}[Bounded holomorphic functions]
\label{coro:bdd-holo}
Let $\hol^\infty(\teich_{g,m})$ be the complex Banach space of bounded holomorphic functions on $\teich_{g,m}$ with the supremum norm. 
Then, the linear mapping
\begin{equation}
\label{eq:isometric_embedding}
\hol^\infty(\teich_{g,m}) \ni f\mapsto f^*\in L^\infty(\partial \Bers{x_0})
\end{equation}
is an isometric embedding.
\end{corollary}

Indeed, the injectivity follows from Theorem \ref{thm:F-M-Riesz}, and the isometricity is deduced from the maximum principle.

Since the Teichm\"uller space is the deformation space of marked Riemann surfaces, the boundary consists of topological data which record how Riemann surfaces degenerate. Hence, the researches with the boundaries of the Teichm\"uller space are expected to contribute to the study of the low-dimensional topology. Indeed, Theorem \ref{thm:F-M-Riesz} deduces the following.

\begin{corollary}[Non-ergodicity of the action the Torelli group on $\mathcal{PMF}$]
\label{coro:2}
The action of the Torelli group $\mathcal{I}_g$ on the space $\mathcal{PMF}=\mathcal{PMF}(\Sigma_g)$ of projective measured foliations on $\Sigma_g$ is not ergodic.
\end{corollary}

A measure on a measure space is said to be \emph{quasi-invariant} under a group action if every element of the group preserves the null sets. A group action on a probability space (whose probability measure is quasi-invariant under the action) is said to be \emph{ergodic} if any invariant measurable set under the action is either null or full (cf. \cite[\S1.5, \S10.6]{Walters:1982}). It is known that the Thurston measure is quasi-invariant with respect to the action of the mapping class group, and that the action of the mapping class group on $\mathcal{PMF}$ is ergodic (cf. \cite{Masur:1982}. See also \cite{Rees:1981}).

Recall that the \emph{Torelli group} $\mathcal{I}_g$ is a subgroup of the mapping class group of a closed surface $\Sigma_g$ of genus $g$ which consists of mapping classes acting trivially on the first homology group on $H_1(\Sigma_g)$ (e.g. \cite{Johnson:1983}).  The Torelli group is known to be a fascinating big subgroup unless $g=1$ ($\mathcal{I}_1$ is trivial). When $g=2$, $\mathcal{I}_2$ is an infinite rank free group, but  $\mathcal{I}_g$ is known to be a finitely generated torsion free group for $g\ge 3$ (cf. \cite{Johnson:1983}, \cite{McCullough_Miller:1986} and \cite{Mess:1992}). Moreover, in contrast with Corollary \ref{coro:2}, the ergodicity for the natural actions of $\mathcal{I}_g$ on the representation spaces are observed in many cases (e.g \cite{Goldman_Xia:2012}).

\section{Teichm\"uller theory}

For the Teichm\"uller theory, see \cite{Abikoff:1980, Imayoshi_Taniguchi:1992, Nag:1988} for instance.

\subsection{Teichm\"uller space}
A \emph{marked Riemann surface} $(M,f)$ of type $(g,m)$
is a pair of a Riemann surface $M$ of analytically finite type $(g,m)$
and an orientation preserving homeomorphism
$f\colon \Sigma_{g,m}\to M$.
Two marked Riemann surfaces $(M_{1},f_{1})$ and $(M_{2},f_{2})$
of type $(g,m)$ are
\emph{(Teichm\"uller) equivalent} if there is a conformal mapping
$h\colon M_{1}\to M_{2}$ such that $h\circ f_{1}$ is homotopic to $f_{2}$.
The \emph{Teichm\"uller space} $\teich_{g,m}$ of type $(g,m)$
is the set of all Teichm\"uller equivalence classes
of marked Riemann surfaces of type $(g,m)$.

\subsection{Teichm\"uller distance}
\label{subsec:teich_distance}
The \emph{Teichm\"uller distance} $d_T$ is a complete distance
on $\teich_{g,m}$ defined by
$$
d_T(x_1,x_2)=\frac{1}{2}\log \inf_hK(h)
$$
for $x_i=(M_i,f_i)$ ($i=1,2$),
where the infimum runs over all quasiconformal mapping
$h\colon M_1\to M_2$ homotopic to $f_2\circ f_1^{-1}$
and $K(h)$ is the maximal dilatation of a quasiconformal mapping $h$.
%

For $x=(M,f)\in \teich_{g,m}$,
we denote by $\mathcal{Q}_{x}$ the complex Banach space
of holomorphic quadratic differentials $q=q(z)dz^{2}$
on $M$ with
$$
\|q\|=\int_{M}|q(z)|\frac{\sqrt{-1}}{2}dz\wedge d\overline{z}
<\infty.
$$
Let $\mathcal{Q}^1_{x}=\{q\in \mathcal{Q}_{x}\mid \|q\|=1\}$ be the unit sphere. 

Let $x=(M,f)\in \teich_{g,m}$. For $q\in \mathcal{Q}_x-\{0\}$ and $t\in [0,\infty)$, let $f_t$ be the quasiconformal mapping on $M_0$ by the Beltrami differential $\tanh(t)\overline{q}/|q|$. We define the \emph{Teichm\"uller (geodesic) ray $\ray_q\colon [0,\infty)\to \teich_{g,m}$ associated to $q$} by $\ray_q(t)=(f_t(M),f_t\circ f)$. Teichm\"uller ray is a geodesic ray with respect to $d_T$. Namely, for $t_1$, $t_2\in [0,\infty)$,
$$
d_T(\ray_q(t_1),\ray_q(t_2))=|t_1-t_2|.
$$

\subsection{Measured foliations and laminations}
Let $\mathcal{S}$ be the set of homotopy classes of 
essential simple closed curves on $\Sigma_{g,m}$.
Let $i(\alpha,\beta)$ denote the \emph{geometric intersection number}
for simple closed curves $\alpha,\beta\in \mathcal{S}$.
Let $\mathcal{WS}=\{t\alpha\mid t\ge 0, \alpha\in \mathcal{S}\}$
be the set of weighted simple closed curves.
The intersection number on $\mathcal{WS}$
is defined by
\begin{equation}
\label{eq:intersection-number-WS}
i(t\alpha,s\beta)=ts\,i(\alpha,\beta)\quad
(t\alpha,s\beta\in \mathcal{WS}).
\end{equation}

\subsubsection{Measured foliations}
We consider an embedding
\begin{equation*}
\mathcal{WS}\ni t\alpha\mapsto [\mathcal{S}\ni \beta
\mapsto i(t\alpha,\beta)]\in \mathbb{R}_{\ge 0}^{\mathcal{S}}.
\end{equation*}
We topologize the function space
$\mathbb{R}_{\ge 0}^{\mathcal{S}}$
with the topology of pointwise convergence.
The closure $\mathcal{MF}$ of the image of the embedding is called
the \emph{space of measured foliations} on $\Sigma_{g,m}$.
Let
$$
{\rm proj}\colon \mathbb{R}_{\ge 0}^{\mathcal{S}}-\{0\}\to \mathbb{PR}_{\ge 0}^{\mathcal{S}}=(\mathbb{R}_{\ge 0}^{\mathcal{S}}-\{0\})/\mathbb{R}_{>0}
$$
be the projection.
The image $\mathcal{PMF}={\rm proj}(\mathcal{MF}-\{0\})$
is called the space of \emph{projective measured foliations}  on $\Sigma_{g,m}$.
We write $[F]$ the projective class of $F\in \mathcal{MF}-\{0\}$.
It is known that $\mathcal{MF}$ and $\mathcal{PMF}$ are homeomorphic 
to $\mathbb{R}^{6g-6+2m}$ and $\mathbb{S}^{6g-7+2m}$,
respectively (cf. \cite{Douady_Fathi_Fried_Laudenbach_Poenaru_Shub:1979}).
By definition,
$\mathcal{MF}$
contains $\mathcal{WS}$
as a dense subset.
The intersection number extends continuously 
as a non-negative function $i(\,\cdot\,,\,\cdot\,)$
on $\mathcal{MF}\times \mathcal{MF}$
satisfying $i(F,F)=0$ and $F(\alpha)=i(F,\alpha)$
for $F\in \mathcal{MF}\subset \mathbb{R}_{\ge 0}^{\mathcal{S}}$
and $\alpha\in \mathcal{S}$.
%

\subsubsection{Measured laminations}
Fix a hyperbolic structure of finite area on $\Sigma_{g,m}$.
A \emph{geodesic lamination} $L$ on $\Sigma_{g,m}$ is a non-empty closed set
which is a disjoint union of complete simple geodesics,
where a geodesic is said to be \emph{complete}
if it is either closed or has infinite length in both of its ends.
The geodesics in $L$ are called the \emph{leaves} of $L$.
A \emph{transverse measure} for a geodesic lamination $L$
means an assignment a Borel measure to each arc transverse to $L$,
subject to the following two conditions:
If the arc $k'$ is contained in the transverse arc $k$,
the measure assigned to $k'$ is the restriction of the measure assigned to $k$;
and
if the two arcs $k$ and $k'$ are homotopic through a family of arcs transverse to $L$,
the homotopy sends the measure assigned to $k$ to the measure assigned to $k'$.
A transverse measure to a geodesic lamination $L$
is said to have \emph{full support} if
the support of the measure assigned to each transverse arc $k$ is exactly $k\cap L$.
A \emph{measured lamination} $L$ is a pair  consisting of a geodesic lamination
called the \emph{support} of $L$, 
and full support transverse measures to the support. 
Let $\mathcal{ML}$ be the set of measured laminations on $\Sigma_{g,m}$
(with fixing a complete hyperbolic structure).
A weighted simple closed curve $t\alpha$ is identified with
a measured lamination consisting of a simple closed geodesic homotopic to $\alpha$
and an assignment $t$-times the Dirac measures
whose support consists of the intersection to transverse arcs.
The intersection number \eqref{eq:intersection-number-WS}
on weighted simple closed curves
extends continuously to $\mathcal{ML}\times \mathcal{ML}$.

It is known that there is a canonical identification
$\mathcal{MF}\cong \mathcal{ML}$ such that
$F\in \mathcal{MF}$ corresponds to $L$ if and only if
\begin{equation*}
i(F,\alpha)=i(L,\alpha)\quad (\alpha\in \mathcal{S})
\end{equation*}
(e.g. \cite{Bonahon:2001} , \cite{Penner_Harer:1992} and \cite{Thurston:1980}).

\begin{quote}{{\bf Convention}}
Henceforth,
we will frequency use the canonical correspondence between measured laminations
and measured foliations.
\end{quote}

For $F\in \mathcal{MF}$,
we denote by $L(F)$ the support of the corresponding measured lamination.
For simplicity,
we call $L(F)$ the \emph{support lamination} of $F$.

An $F\in \mathcal{MF}$ is called \emph{minimal}
if any leaf of $L(F)$ is dense in $L(F)$
(with respect to the induced topology from $\Sigma_{g,m}$).
An $F\in \mathcal{MF}$ is called \emph{filling}
if any complementary region of $L(F)$ is either an ideal polygon or 
a once punctured ideal polygon, which is equivalent to say that $i(F,\alpha)\ne 0$ for all $\alpha\in \mathcal{S}$ (e.g. \cite[\S2.2]{Lindenstrauss_Mirzakhani:2008}).
A measured lamination $L$ is said to be \emph{uniquely ergodic} 
if $L'\in \mathcal{ML}$ satisfies $i(L,L')=0$,
then $L'=tL$ for some $t\ge 0$.
A measured foliation is said to be \emph{uniquely ergodic} if
so is the corresponding measured lamination.
%

\subsection{Hubbard Masur differentials and extremal length}
Let $x=(M,f)\in \teich_{g,m}$ and $q\in \mathcal{Q}_x$. We can define the \emph{vertical foliation} $v(q)\in \mathcal{MF}$ of $q$ by
\begin{equation}
\label{eq:vertical_foliation}
i(v(q),\alpha)=\inf_{\alpha'\sim f(\alpha)}\int_{\alpha'}|{\rm Re}(\sqrt{q})|
\end{equation}
Hubbard and Masur \cite{Hubbard_Masur:1979} showed that for $x=(M,f)\in \teich_{g,m}$ and $F\in \mathcal{MF}$, there is a unique $q_{F,x}\in \mathcal{Q}_x$ such that $v(q_{F,x})=F$.
In fact, for $x\in \teich_{g,m}$, the correspondence $\mathcal{MF}\ni F\mapsto q_{F,x}\in \mathcal{Q}_x$ is homeomorphic.
We call the differential $q_{F,x}$ the \emph{Hubbard-Masur differential} for $F$ on $x$. 

For $F\in \mathcal{MF}$,
we define the \emph{extremal length} of $F$ on $x=(M,f)\in \teich_{g,m}$ by
$$
\ext_x(F)=\|q_{F,x}\|=\int_M |q_{F,x}(z)|dxdy.
$$

\subsection{Bers slice}
\label{subsec:Bers-slice}
Fix $x_{0}=(M_{0},f_{0})\in \teich_{g,m}$
and let $\Gamma_{0}$ be the marked Fuchsian group acting on $\mathbb{H}$
uniformizing $M_{0}$ with the marking $\pi_1(\Sigma_{g,m})\cong \Gamma_0$
induced by $f_0$.
Let $A_{2}(\mathbb{H}^{*},\Gamma_{0})$ be the Banach space of automorphic forms
on $\mathbb{H}^{*}=\hat{\mathbb{C}}-\overline{\mathbb{H}}$
of weight $-4$ with the hyperbolic supremum norm.
For each $\varphi\in A_{2}(\mathbb{H}^{*},\Gamma_{0})$,
we can define a locally univalent meromorphic mapping $W_{\varphi}$
on $\mathbb{H}^{*}$ and the monodromy homomorphism 
$\rho_{\varphi}\colon \Gamma_{0}\to {\rm PSL}_{2}(\mathbb{C})$
such that 
the Schwarzian derivative of
$W_{\varphi}$ 
is equal to $\varphi$ and $\rho_{\varphi}(\gamma)\circ W_{\varphi}=W_{\varphi}\circ \gamma$
for all $\gamma\in \Gamma_{0}$.
Let $\Gamma_{\varphi}=\rho_{\varphi}(\Gamma_{0})$. Notice that all group $\Gamma_\varphi$ is \emph{marked} with a surjective homomorphism $\rho_\varphi\colon \Gamma_0(\cong \pi_1(\Sigma_{g,m}))\to \Gamma_\varphi$.

The \emph{Bers slice} $\Bers{x_{0}}$ with base point
$x_{0}\in \teich_{g,m}$
is a domain in $A_{2}(\mathbb{H}^{*},\Gamma_{0})$
which consists of $\varphi\in A_{2}(\mathbb{H}^{*},\Gamma_{0})$
such that $W_{\varphi}$ admits a quasiconformal extension
to $\hat{\mathbb{C}}$.
The Bers slice $\Bers{x_{0}}$ is bounded and
identified biholomorphically
with $\teich_{g,m}$.
Indeed,
any $x\in \teich_{g,m}$ corresponds to
$\varphi$ such that $\Gamma_{\varphi}$ is
the marked quasifuchsian group uniformizing $x_{0}$ and $x$
(cf. \cite{Bers:1960}).
The closure $\overline{\Bers{x_{0}}}$ of $\Bers{x_{0}}$ in $A_{2}(\mathbb{H}^{*},\Gamma_{0})$
is called the \emph{Bers compactification}
of $\teich_{g,m}$.
The boundary $\partial\Bers{x_{0}}$ is called the \emph{Bers boundary}.
For $\varphi\in \overline{\Bers{x_{0}}}$,
$\Gamma_\varphi$ is a marked Kleinian surface group with isomorphism
$\rho_\varphi\colon \pi_1(\Sigma_{g,m})\cong \Gamma_0\to \Gamma_\varphi$.

\subsection{Boundary groups without APTs}
\label{subsec:BoundarygroupswithoutAPT}
A boundary point $\varphi\in \partial\Bers{x_0}$ is called a \emph{cusp} if
there is a non-parabolic element $\gamma\in \Gamma_{0}$
such that $\rho_{\varphi}(\gamma)$ is parabolic
(cf. \cite{Bers:1970}).
Such $\gamma$ or $\rho_{\varphi}(\gamma)$
is called an \emph{accidental parabolic transformation} (APT)
of $\varphi$ or $\Gamma_{\varphi}$.
Let $\partial^{cusp}\Bers{x_{0}}$ be the set of cusps in $\partial\Bers{x_{0}}$
and set
$\partial^{mf}\Bers{x_{0}}
=\partial\Bers{x_{0}}-\partial^{cusp}\Bers{x_{0}}$.

For $\varphi\in \partial^{mf}\Bers{x_{0}}$,
the quotient manifold $\mathbb{H}^{3}/\Gamma_{\varphi}$
has two (non-cuspidal) ends
corresponding to $\Sigma_{g,m}\times (0,\infty)$
and $\Sigma_{g,m}\times (-\infty,0)$.
The negative end is geometrically finite and the surface at infinity is conformally equivalent to the mirror of $M_{0}$.
To another end, we assign a unique minimal and filling geodesic lamination,
called the \emph{ending lamination} for $\varphi$
 (cf. \cite{Bonahon:1986} and \cite{Thurston:1980}).

Let $x_0\in \teich_{g,m}$.
Let $\mathcal{PMF}^{mf}$ be the set of projective classes
of minimal and filling measured foliations.
By virtue of the Ending Lamination Theorem and the Thurston double limit theorem,
we have the closed continuous surjective mapping
\begin{equation}
\label{eq:curve-complex-mininal-bers}
\Xi_{x_0}\colon \mathcal{PMF}^{mf}\to 
\partial^{mf} \Bers{x_0}
\end{equation}
which assigns $[F]\in \mathcal{PMF}^{mf}$
to the boundary group whose ending lamination is equal to $L(F)$ (cf. \cite{Brock_Canary_Minsky:2012}).
The preimage of any point in $\partial^{mf} \Bers{x_0}$ is compact
(cf. \cite{Leininger_Schleimer:2009}).
$\mathcal{PMF}^{mf}$ contains a subset $\mathcal{PMF}^{ue}$ consisting of minimal, filling and uniquely ergodic measured foliations.
Let $\partial^{ue} \Bers{x_0}$ be the image of $\mathcal{PMF}^{ue}$ under the identification \eqref{eq:curve-complex-mininal-bers}.

\subsection{Teichm\"uller rays associated to projective measured foliations}
For $[F]\in \mathcal{PMF}$ and $x\in \teich_{g,m}$,
let $\ray^x_F\colon [0,\infty)\to \teich_{g,m}$ be the Teichm\"uller ray associated to $q_{F,x}$.  Namely, $\ray^x_F=\ray_{q_{F,x}}$. The ray $\ray^x_F$ is independent of the choice of the representative in the class $[F]$ (cf. \S\ref{subsec:teich_distance}).

The following proposition folllows from the ending lamination theorem \cite{Brock_Canary_Minsky:2012} and the continuity of the length of laminations (cf. \cite{Brock:2000}, \cite{Ohshika:1990}. See also \cite[Theorem 6.1]{Brock:2001}).

\begin{proposition}
\label{prop:Teichmuller-limit}
Let $x_0\in \teich_{g,m}$.
For $x\in \teich_{g,m}$ and  $[H]\in \mathcal{PMF}^{mf}$,
the Teichm\"uller ray $\ray_{H}^x$ converges to the totally degenerate
group $\varphi_H$ without APT in $\partial^{mf}\Bers{x_{0}}$
whose ending lamination is $L(H)$.
\end{proposition}

\begin{proof}
We give a brief proof of Proposition \ref{prop:Teichmuller-limit} for reader's convenience.
From the (analytic) definition of the extremal length, the Bers inequality (\cite[Theorem 3]{Bers:1970}) and the continuity of the hyperbolic length of the measured foliations (\cite[Theorem 2]{Brock:2000}), 
we have $\ell_{\varphi_t}(H)\le C\ext_{x_t}(H)^{1/2}\le Ce^{-t}\ext_{x}(H)^{1/2}$
for $t\ge 0$, where $C$ is a positive constant depending only on $g$ and $m$, $\varphi_t\in \Bers{x_0}$ corresponds to $\ray^{x}_H(t)\in \teich_{g,m}$, and $\ell_{\varphi_t}(H)$ is the hyperbolic length of $H$ for the (marked) quasi-Fuchsian manifold $\mathbb{H}^3/\Gamma_{\varphi_t}$. Let $\varphi\in \partial\Bers{x_0}$ be an accumulation point of the Teichm\"uller ray $\ray^x_H$. Letting $t\to \infty$, we get $\ell_{\varphi}(H)=0$, which implies that $H$ is non-realizable in $\mathbb{H}^3/\Gamma_\varphi$. Since $H$ is minimal and filling, the support $L(H)$ is the ending lamination of $\mathbb{H}^3/\Gamma_\varphi$. Since the Bers slice $\Bers{x_0}$ is the deformation space of (marked) quasi-Fuchsian manifolds with fixing one end to be $x_0$ (cf. \cite{Bers:1960}), from the ending lamination theorem, we conclude $\varphi=\varphi_H$.
\end{proof}

Let $[\omega]$ be a mapping class on $\Sigma_g$ and $[\omega]_*$ is the action on $\teich_{g,m}$ induced by $[\omega]$ (e.g. \cite{Imayoshi_Taniguchi:1992}).
Then, from \eqref{eq:vertical_foliation}
$$
[\omega]_*\circ \ray^{x}_F(t)=\ray^{[\omega]_*(x)}_{\omega(F)}(t)
$$
for $t\ge 0$.
Since $[\omega]_*$ naturally extends to $\partial^{mf}\Bers{x_0}$
(cf. \cite{Bers:1981} and \cite{Brock_Canary_Minsky:2012}),
\begin{equation}
\label{eq:mapping_class_ray}
\lim_{t\to \infty}[\omega]_*\circ \ray^{x}_F(t)=\varphi_{\omega(F)}=[\omega]_*(\varphi_F).
\end{equation}

\subsection{Thurston measure}
\label{subsec:Thurston_measure}
There is a unique (up to constant multiple) locally finite mapping class group-invariant ergodic measure $\ThursM$ on $\mathcal{MF}$ supported on the sets of filling measured foliations. The measure $\ThursM$ is called the \emph{Thurston measure} (cf. \cite[Theorem 7.1]{Lindenstrauss_Mirzakhani:2008}). 
For $x\in \teich_{g,m}$ and $E\subset \mathcal{PMF}$,
we set
$$
{\rm Cone}(E)_x=\left\{t\dfrac{F}{\ext_x(F)^{1/2}}\in \mathcal{MF}\mid [F]\in E, \ 0\le t\le 1\right\}.
$$
We define a probability measure $\ThursM^x$ on $\mathcal{PMF}$ by
$$
\ThursM^x(E)=\dfrac{\ThursM({\rm Cone}(E)_x)}{\ThursM({\rm Cone}(\mathcal{PMF})_x)}
$$
for $E\subset \mathcal{PMF}$. For simplicity, we also call $\ThursM^x$ the \emph{Thurston measure} on $\mathcal{PMF}$ associated to $x\in \teich_{g,m}$.

\section{Radial limit theorem}
\label{sec:radial_limit_theorem}
In this section, we shall show the following:

\begin{theorem}[Radial limit theorem]
\label{thm:radial_limit_theorem}
For a bounded pluriharmonic function $u$ on $\teich_{g,m}$, there is a full-measure set $\mathcal{E}_0=\mathcal{E}_0(u)\subset \mathcal{PMF}$ depending only on $u$ with respect to the Thurston measure with the following properties:
\begin{itemize}
\item[{\rm (1)}]
each element in $\mathcal{E}_0$ is minimal, filling and uniquely ergodic;
\item[{\rm (2)}]
the radial limit $\lim_{t\to \infty}u(\ray^x_F(t))$ exists for all $x\in \teich_{g,m}$ and $[F]\in \mathcal{E}_0$; and
\item[{\rm (3)}]
the radial limit is independent of the choice of the base point. Namely,
$$
\lim_{t\to \infty}u(\ray^{x_1}_F(t))=\lim_{t\to \infty}u(\ray^{x_2}_F(t))
$$
for $[F]\in \mathcal{E}_0$ and $x_1,x_2\in \teich_{g,m}$.
\end{itemize}
\end{theorem}

Following Theorem \ref{thm:radial_limit_theorem}, we define a bounded measurable function on $\partial\Bers{x_0}$ by
\begin{equation}
\label{eq:radial_limit}
u^*(\varphi_F)=\begin{cases}
{\displaystyle \lim_{t\to \infty}u(\ray^{x_0}_F(t))} & ([F]\in \mathcal{E}_0) \\
0 & ([F]\in \mathcal{PMF}\setminus \mathcal{E}_0)
\end{cases}
\end{equation}
for a bounded pluriharmonic function $u$ on $\teich_{g,m}$, where $\varphi_F\in \partial \Bers{x_0}$ is the boundary group with ending lamination $L(F)$, $\mathcal{E}_0$ is a full measure set in $\mathcal{PMF}$ with respect to the Thurston measure defined in Theorem \ref{thm:radial_limit_theorem} for $u$. Since $\lim_{t\to \infty}\ray^{x}_F(t)=\varphi_F$ for all $x\in \teich_{g,m}$ (cf. Proposition \ref{prop:Teichmuller-limit}), the \emph{radial limit} $u^*$ of $u$ is independent of the choice of $x_0\in \teich_{g,m}$.
In particular,
from \eqref{eq:mapping_class_ray}, the radial limit is \emph{natural} with respect to the action of the mapping class group in the sense that
\begin{equation}
\label{eq:naturality_MCG}
(u\circ [\omega]_*)^*=u^*\circ [\omega]_*
\end{equation}
for a mapping class $[\omega]$ on $\Sigma_g$.

\subsection{Projectification of $\mathcal{MF}$ and Disintegration}
\label{subsec:projectification_PMF}
Fix $x_0\in \teich_{g,m}$. Let $\mathbb{S}^1=\{|z|=1\}$ be the unit circle.
We define the action of $\mathbb{S}^1$ on $\mathcal{PMF}$ by
$$
\mathbb{S}^1\times \mathcal{PMF}\ni (e^{i\alpha},[F])\mapsto A_\alpha([F]):=[v(e^{i\alpha}q_{F,x_0})]\in \mathcal{PMF}.
$$
We denote by $\mathbb{P}_{x_0}\mathcal{MF}$ the hopf quotient $\mathcal{PMF}/\mathbb{S}^1\cong \mathbb{S}^{6g-7+2m}/\mathbb{S}^1\cong \mathbb{C}\mathbb{P}^{3g-4+m}$ and by $\Pi^{x_0}$ the projection $\mathcal{PMF}\to \mathbb{P}_{x_0}\mathcal{MF}$. Let $\nu^{x_0}$ be the push forward measure of $\mu_{Th}^{x_0}$ via the projection. By definition, $\nu^{x_0}$ is a probability measure on $\mathbb{P}_{x_0}\mathcal{MF}$.

From the disintegration theorem, there is the disintegration $\{\lambda_{t}\mid t\in \mathbb{P}_{x_0}\mathcal{MF}\}$ with respect to the projection (cf. \cite[Theorem 1]{Chang_Polland:1997}). Namely, 
each $\lambda_t$ is a finite measure on $\mathbb{P}_{x_0}\mathcal{MF}$ concentrated on $(\Pi^{x_0})^{-1}(t)$ (i.e. $\lambda_t(\{[F]\in \mathcal{PMF}\mid \Pi^{x_0}([F])\ne t\})=0$); for each nonnegative measurable function $f$ on $\mathcal{PMF}$,
\begin{itemize}
\item[(i)]
$\displaystyle\mathbb{P}_{x_0}\mathcal{MF}\ni t\mapsto \int_{\mathcal{PMF}}f([F])\,d\lambda_t([F])$ is measurable;
\item[(ii)]
$\displaystyle \int_{\mathbb{P}_{x_0}\mathcal{MF}}\left(\int_{\mathcal{PMF}}f([F])\,d\lambda_t([F])\right)d\nu^{x_0}(t)=\int_{\mathcal{Q}^1_{x_0}}f([F])\,d\ThursM^{x_0}([F])$.
\end{itemize}
Furthermore, the measures $\{\lambda_t\}_t$ are determined up to an almost sure equivalence in the sense that if $\{\lambda_t^{*}\}_t$ is another disintegration, then $\nu^{x_0}(\{t\in \mathbb{P}_{x_0}\mathcal{MF}\mid \lambda^*_t\ne \lambda_t\})=0$.
From \cite[Theorem 2]{Chang_Polland:1997}, $\lambda_t$ is a probability measure for almost all $t\in \mathbb{P}_{x_0}\mathcal{MF}$.

By definition,
for any $t\in \mathbb{P}_{x_0}\mathcal{MF}$, there is a canonical identification
\begin{equation}
\label{eq:identification_Pi}
\mathbb{S}^1\ni e^{i\theta}\mapsto [v(e^{i\theta}q_{F,x_0})]\in (\Pi^{x_0})^{-1}(t)
\end{equation}
for all $[F]\in \mathcal{PMF}$ with $t=\Pi^{x_0}([F])$. The identification \eqref{eq:identification_Pi} is determined up to composing rotations on $\mathbb{S}^1$. For $t\in \mathbb{P}_{x_0}\mathcal{MF}$, we denote by $\Theta_t$ the push-forward measure of $d\theta/2\pi$ on $(\Pi^{x_0})^{-1}(t)$ via the identification \eqref{eq:identification_Pi}. Since the measure $d\theta/2\pi$ on $\mathbb{S}^1$ is invariant under the rotation on $\mathbb{S}^1$, the measure $\Theta_t$ is well-defined independently of  the choice of $q$ in the identification \eqref{eq:identification_Pi}.

\begin{proposition}
\label{prop:absolutecontiuous}
 For almost all $t\in \mathbb{P}_{x_0}\mathcal{MF}$,
 $\lambda_t=\Theta_t$.
\end{proposition}

\begin{proof}
Fix $\alpha\in [0,2\pi)$.
Dumas \cite[Corollary 5.9]{Dumas:2015} shows that the action $A_\alpha$ preserves the Thurston measure.
Namely, $(A_\alpha)_*\ThursM^{x_0}=\ThursM^{x_0}$ fir all $x\in \teich_{g,m}$ (Dumas treated the case where $m=0$ and $g\ge 2$. However the proof is also available for $m>0$ with $2g-2+m>0$).
Hence, for a non-negative measurable function $f$ on $\mathcal{PMF}$,
$$
\mathbb{P}_{x_0}\mathcal{MF}\ni t\mapsto
\int_{\mathcal{PMF}}f([F])\,d((A_\alpha)_*\lambda_t)([F])
=
\int_{\mathcal{PMF}}f\circ A_\alpha([F])\,d\lambda_t([F])
$$
is measurable
and
\begin{align*}
\int_{\mathcal{PMF}}f([F])\,d\ThursM^{x_0}([F])
&=
\int_{\mathcal{PMF}}f([F])\,d((A_\alpha)_*\ThursM^{x_0})([F]) \\
&=
\int_{\mathcal{PMF}}f\circ A_\alpha([F])\,d\ThursM^{x_0}([F])
\\
&= \int_{\mathbb{P}_{x_0}\mathcal{MF}}\left(\int_{\mathcal{PMF}}f\circ A_\alpha(q)\,d\lambda_t([F])\right)d\nu^{x_0}(t) \\
&= \int_{\mathbb{P}_{x_0}\mathcal{MF}}\left(
\int_{\mathcal{PMF}}f([F])\,d((A_\alpha)_*\lambda_t)([F])
\right)d\nu^{x_0}(t)
\end{align*}
from the property (ii) of the disintegration discussed above.
Therefore, $\{(A_\alpha)_*\lambda_t\}_t$ is also the disintegration with respect to the projection $\Pi^{x_0}$ for all $\alpha\in [0,2\pi)$. 

Now, we assume that $\alpha/2\pi$ is irrational. From the uniqueness of the disintegration, $(A_\alpha)_*\lambda_t=\lambda_t$ almost everywhere on $\mathbb{P}_{x_0}\mathcal{MF}$. This means that $\lambda_t$ is an invariant measure on $\mathbb{S}^1$ in terms of the irrational rotation $A_\alpha$. Since any irrational rotation has no periodic points in $\mathbb{S}^1$, the rotation $A_\alpha$ is uniquely ergodic  (cf. \cite[Theorem 6.18]{Walters:1982}). Hence, the invariant measure $\lambda_t$ coincides with a constant multiple of the Lebesgue measure. Since $\lambda_t$ is a probability measure, we conclude that $\lambda_t=\Theta_t$ almost all $t\in \mathbb{P}_{x_0}\mathcal{MF}$.
\end{proof}

\subsection{Proof of Theorem \ref{thm:radial_limit_theorem}}
\label{subsec:proof_radial_limit_theorem}
Let $u$ be a bounded pluriharmonic function on $\teich_{g,m}$.
For $n,m\in \mathbb{N}$, we define
\begin{align*}
\mathcal{E}_{n,m}
&=
\left\{[F]\in \mathcal{PMF}\mid
\limsup_{k\to \infty,k\in \mathbb{N}}u(\ray^{x_0}_F(k/2^m))
\le
\liminf_{k\to \infty,k\in \mathbb{N}}u(\ray^{x_0}_F(k/2^m))
+
\dfrac{1}{n}
\right\}
\\
\mathcal{E}_\infty
&=\cap_{n>0}\left(\cap_{m>0}\mathcal{E}_{n,m}\right).
\end{align*}
We notice that $\mathcal{E}_{n_2,m_2}\subset \mathcal{E}_{n_1,m_1}$
for $n_2\ge n_1$ and $m_2\ge m_1$. Indeed, for $[F]\in \mathcal{E}_{n_2,m_2}$,
\begin{align*}
\limsup_{k\to \infty}u(\ray^{x_0}_F(k/2^{m_1}))
&\le \limsup_{k\to \infty}u(\ray^{x_0}_F(k/2^{m_2}))
\le\liminf_{k\to \infty}u(\ray^{x_0}_F(k/2^{m_2}))+\dfrac{1}{n_2} \\
&\le \liminf_{k\to \infty}u(\ray^{x_0}_F(k/2^{m_1}))+\dfrac{1}{n_1},
\end{align*}
and hence $[F]\in \mathcal{E}_{n_1,m_1}$.
Since $\mathcal{PMF}\ni [F]\mapsto u(\ray^{x_0}_F(k/2^m))$ is continuous for fixed $k$ and $m$, each $\mathcal{E}_{n,m}$ is measurable. Hence, $\mathcal{E}_\infty$ is also measurable.

We claim 
\begin{lemma}
\label{lemma:1}
For $[F]\in \mathcal{PMF}$,
$[F]\in \mathcal{E}_\infty$ if and only if the limit $\displaystyle\lim_{t\to \infty}u(\ray^{x_0}_F(t))$ exists. 
\end{lemma}

\begin{proof}[Proof of Lemma \ref{lemma:1}]
Suppose that $[F]\in \mathcal{E}_\infty$.
Let $n,m\in \mathbb{N}$. Since $[F]\in \mathcal{E}_{n,m}$, from the Schwarz lemma discussed in \S\ref{subsec:Schwarz_lemma} below, for any $t>0$, there is $k\in \mathbb{N}$ such that
$$
|u(\ray^{x_0}_F(t))-u(\ray^{x_0}_F(k/2^{m}))|\le \dfrac{C}{2^m},
$$
where $C>0$ is a constant depending only on $\|u\|_\infty$.
Therefore, we get
\begin{align*}
\limsup_{t\to \infty}u(\ray^{x_0}_F(t))&\le \limsup_{k\to \infty}u(\ray^{x_0}_F(k/2^{m}))+\dfrac{C}{2^m} \\
&\le\liminf_{k\to \infty}u(\ray^{x_0}_F(k/2^{m}))+\dfrac{C}{2^m}+\dfrac{1}{n} \\
&\le \liminf_{t\to \infty}u(\ray^{x_0}_F(t))+\dfrac{C}{2^{m-1}}+\dfrac{1}{n}
\end{align*}
for all $m$. Since $[F]\in \cap_{m>0}\mathcal{E}_{n,m}$ for all $n$,
by letting $m\to \infty$, we have
$$
\limsup_{t\to \infty}u(\ray^{x_0}_F(t))\le \liminf_{t\to \infty}u(\ray^{x_0}_F(t))+\dfrac{1}{n}.
$$
Since $[F]\in \mathcal{E}_\infty$,
by letting $n\to \infty$, we conclude that the limit of $u$ along the Teichm\"uller ray $\ray^{x_0}_F$ exists.

Conversely, assume that the limit of $u$ along the Teichm\"uller ray $\ray^{x_0}_F$ exists. Let $n,m\in \mathbb{N}$. Then
\begin{align*}
\limsup_{k\to \infty}u(\ray^{x_0}_F(k/2^{m}))
&\le \limsup_{t\to \infty}u(\ray^{x_0}_F(t))
=\liminf_{t\to \infty}u(\ray^{x_0}_F(t)) \\
&\le \liminf_{k\to \infty}u(\ray^{x_0}_F(k/2^{m}))\le \liminf_{k\to \infty}u(\ray^{x_0}_F(k/2^{m}))+\dfrac{1}{n}
\end{align*}
and hence $[F]\in \mathcal{E}_{n,m}$. Threfore $[F]\in \mathcal{E}_\infty$.
\end{proof}

Let us return to the proof of Theorem \ref{thm:radial_limit_theorem}.
Let $x\in \teich_{g,m}$ and $[F]\in \mathcal{PMF}$. Consider the Teichm\"uller disk $\mathcal{R}_{[F]}\colon \mathbb{D}\to \teich_{g,m}$ associated to $q_{F,x_0}$, which does not depend on the choice of representatives in the projective class $[F]$. Since $u\circ \mathcal{R}_{[F]}$ is a bounded harmonic function on $\mathbb{D}$, there is a full measure set $E_{[F]}\subset \mathbb{S}^1=\partial \mathbb{D}$ with respect to $d\theta/2\pi$ such that the radial limit of $u$ exists along the Teichm\"uller ray defined by $e^{i\theta}q_{F,x_0}$ and $e^{i\theta}\in E_{[F]}$ by Fatou's theorem (cf. \cite[Theorem VI.6, Chapter IV]{Tsuji:1959}). 

For $t\in \mathbb{P}_{x_0}\mathcal{MF}$, we fix $[F_t]\in \mathcal{PMF}$ with $\Pi^{x_0}([F_t])=t$. Notice for the later argument that $E_{[H_1]}=e^{i\alpha}E_{[H_2]}$ for some $\alpha\in \mathbb{R}$ when $\Pi^{x_0}([H_1])=\Pi^{x_0}([H_2])$.
Let $\charfunc_{\mathcal{E}_\infty}$ be the characteristic function of $\mathcal{E}_\infty$ on $\mathcal{PMF}$.
By Proposition \ref{prop:absolutecontiuous}, $E_{[F_t]}$ is a full measure set in $(\Pi^{x_0})^{-1}(t)$ with respect to $\lambda_t$ for almost all $t\in \mathbb{P}_{x_0}\mathcal{MF}$. From Lemma \ref{lemma:1},
$$
\int_{\mathcal{PMF}}\charfunc_{\mathcal{E}_\infty}(q)\,d\lambda_t(q)\ge 
\lambda_t(\{[v(e^{i\theta}q_{F_t,x_0})]\mid \theta\in E_{[F_t]}\})
=\int_{\mathcal{PMF}}\,d\lambda_t(q)
$$
for almost all $t\in \mathbb{P}_{x_0}\mathcal{MF}$. From the disintegration theorem, we obtain
\begin{align*}
1\ge \ThursM^{x_0}(\mathcal{E}_\infty)
&=\int_{\mathcal{PMF}}\charfunc_{\mathcal{E}_\infty}(q)d\ThursM^{x_0}(q) \\
&=\int_{\mathbb{P}_{x_0}\mathcal{MF}}\left(\int_{\mathcal{PMF}}\charfunc_{\mathcal{E}_\infty}(q)\,d\lambda_t(q)\right)d\nu^{x_0}(t) \\
&\ge 
\int_{\mathbb{P}_{x_0}\mathcal{MF}}\left(\int_{\mathcal{PMF}}d\lambda_t(q)\right)d\nu^{x_0}(t) 
= \ThursM^{x_0}(\mathcal{PMF})=1.
\end{align*}
This implies that $\mathcal{E}_\infty$ is a full measure set in $\mathcal{PMF}$ with respect to the Thurston measure $\ThursM^{x_0}$. 

We define
$$
\mathcal{E}_0=\{[F]\in \mathcal{E}_{\infty}\mid \mbox{$F$ is minimal, filling and uniquely ergodic}\}.
$$
From \cite[Theorem 2]{Masur:1982} and the above discussion, $\mathcal{E}_0$ is a full measure set in $\mathcal{PMF}$ with respect to the Thurston measure $\ThursM^{x_0}$. Since $\ThursM^x$ is absolutely continuous with respect to $\ThursM^{x_0}$ for all $x\in \teich_{g,m}$, the set $\mathcal{E}_\infty$ is also a full-measure set with respect to $\ThursM^{x}$ for all $x\in \teich_{g,m}$ (cf. \cite[\S2.3.1]{Athreya_Bufetov_Eskin_Mirzakhani:2012}).

Let $[F]\in \mathcal{E}_0$ and $x_1\in \teich_{g,m}$. Take an arbitrary small constant $\epsilon>0$. For $t>0$, we take $s(t)>0$ such that
\begin{equation}
\label{eq:proof_fatou-2}
d_T(\ray_F^{x_1}(t),\ray_F^{x_0}(s(t)))\le \inf_{x\in \ray_F^{x_0}([0,\infty))}d_T(x,\ray_F^{x_1}(t))+\epsilon.
\end{equation}
From the Schwarz lemma discussed in  \S\ref{subsec:Schwarz_lemma} below,

\begin{equation}
\label{eq:proof_fatou-1}
\left|u(\ray_F^{x_1}(t))-u(\ray_F^{x_0}(s(t)))\right|\le 
C\inf_{x\in \ray_F^{x_0}([0,\infty))}d_T(x,\ray_F^{x_1}(t))+C\epsilon,
\end{equation}
where the constant $C>0$ is dependent only on $\|u\|_\infty$.
Since $F$ is filling and uniquely ergodic, by \cite[Theorem 2]{Masur:1980}, the first term of the right-hand side in \eqref{eq:proof_fatou-1} tends to $0$ as $t\to \infty$.
In particular, we also obtain $s(t)\to \infty$ as $t\to \infty$ by \eqref{eq:proof_fatou-2}. Since $[F]\in\mathcal{E}_\infty$, $u(\ray^{x_0}_F(s(t)))$ converges to the radial limit  $\lim_{t\to \infty}u(\ray^{x_0}_F(t))$ as $t\to \infty$ by Lemma \ref{lemma:1}. Therefore, the radial limit $\lim_{t\to \infty}u(\ray^{x_1}_F(t))$ also exists and satisfies
$$
\lim_{t\to \infty}u(\ray^{x_0}_F(t))=\lim_{t\to \infty}u(\ray^{x_1}_F(t))
$$
from \eqref{eq:proof_fatou-1}, since $\epsilon>0$ is taken arbitrary.
This means that $\mathcal{E}_0$ satisifies the properties which we desired.

We finally confirm that the function $u^*$ defined as \eqref{eq:radial_limit} is in $L^\infty(\partial\Bers{x_0})$. Since $u$ is bounded,
so is $u^*$. Hence, we should show that $u^*$ is measurable with respect to the pluriharmonic measure.
Since $\mathcal{E}_0$ is measurable and $\mathcal{PMF}\ni [F]\mapsto u(\ray^{x_0}_F(t))$ is continuous on $\mathcal{PMF}$ for each fixed $t$,
\begin{equation}
\hat{u}^*([F])=
\label{eq:radial_limit2}
\begin{cases}
{\displaystyle \lim_{t\to \infty}u(\ray^{x_0}_F(t))} & ([F]\in \mathcal{E}_0) \\
0 & ([F]\in \mathcal{PMF}\setminus \mathcal{E}_0)
\end{cases}
\end{equation}
is bounded and measurable on $\mathcal{PMF}$ with respect to the Thurston measure $\ThursM^{x_0}$.
Notice that $\Xi_{x_0}$ defined in  \eqref{eq:curve-complex-mininal-bers} is homeomorphic on $\mathcal{E}_0$ onto the image.
Since $u^*\circ \Xi_{x_0}=\hat{u}^*$ and the pushforward $(\Xi_{x_0})_*(\ThursM^{x_0})$ coincides with the pluriharmonic measure on $\partial \Bers{x_0}$,
$u^*$ is a measurable function on $\partial \Bers{x_0}$ with respect to the pluriharmonic measure (cf. \cite[Theorem 1.1]{Miyachi:2023}. See also Demailly \cite{Demailly:1987})).

\subsection{Schwarz lemma}
\label{subsec:Schwarz_lemma}
In the proof of Theorem \ref{thm:radial_limit_theorem}, we use a version of the Schwarz lemma for bounded pluriharmonic functions on a simply connected Kobayashi hyperbolic domain $D\subset \mathbb{C}^n$. The Schwarz lemma discussed here might be well-known. However, we give a brief proof for completeness.

\begin{lemma}[Schwarz lemma]
\label{lem:schwarz_lemma}
Let $u$ be a bounded pluriharmonic function on a simply connected Kobayashi hyperbolic domain $D\subset \mathbb{C}^n$. Then,
$$
|u(z)-u(w)|\le Cd_{D}(z,w)
$$
for $z,w\in D$,
where $d_D$ is the Kobayashi hyperbolic distance on $D$ and $C>0$ is a constant depending only on the sup norm $\|u\|_\infty$ of $u$.
\end{lemma}

 \begin{proof}
 Set $M=\|u\|_\infty$. Since $D$ is simply connected, there is a holomorphic function $f$ on $D$ such that $u={\rm Re}(f)$ (cf. \cite[Theorem 3 in \S K]{Gunning:1990}). In particular $f$ is a holomorphic map from $D$ into a vertical strip $S=\{|{\rm Re}(w)|<M+1\}$. Since vertical translations are conformal automorphisms of $S$, the density of the hyperbolic metric on $S$ at any $w\in S$ is dependent only on the real part ${\rm Re}(w)$. Hence, the vertical projection from $S$ to an open interval $(-M-1,M+1)$ is a contraction with respect to the hyperbolic metric on $S$. By the distance-decreasing property of the Kobayashi metric,
$$
d_S(u(z_1),u(z_2))\le d_S(f(z_1),f(z_2))\le d_{D}(z,w)
$$
for $z,w\in D$, where $d_S$ is the hyperbolic distance on $S$.
Since the image of $u$ is contained in the closed interval $[-M,M]\subset H$, the distance $d_S(u(z_1),u(z_2))$ is comparable with the difference $|u(z_1)-u(z_2)|$ with constants depending only on the bound $M$.
\end{proof}

\section{Identity theorem}
The original F. and M.Riesz theorem is stated as folllows : Let $f$ be a bounded holomorphic function on $\mathbb{D}$.
Suppose that the radial limit (non-tangential limit) $f^*$ of $f$ vanishes on a non-null measurable set in $\partial \mathbb{D}$ with respect to the angle measure. Then, $f$ vanishes (cf. \cite[p.137, Theorem IV.9]{Tsuji:1959}). In this section, we prove Theorem \ref{thm:F-M-Riesz}, which is thought of as a version of F. and M. Riesz theorem for the Teichm\"uller spaces, and Corollary \ref{coro:2}.

\subsection{Proof of Theorem \ref{thm:F-M-Riesz}}
We suppose that there are a non-null measurable set $A\subset \mathcal{PMF}$ and $c\in \mathbb{C}$ such that $f^*\equiv c$ on $\Xi_{x_0}(A)\subset \partial\Bers{x_0}$. We may assume that $A\subset \mathcal{E}_0=\mathcal{E}_0(f)=\mathcal{E}_0({\rm Re}(f))\cap \mathcal{E}_0({\rm Im}(f))$.
By considering $f-c$ instead of $f$, we show only the case where $c=0$.
Furtheremore, since the base point $x_0\in \teich_{g,m}$ of the Bers slice is taken arbitrary in the beginning, it suffices to show that $f(x_0)=0$.
Indeed, let $x_1\in \teich_{g,m}$. Since the radial limit is independent of the choice of the base point, the radial limit of $f$ vanishes on $\Xi_{x_1}(A)\subset \partial\Bers{x_1}$ when we recognize $f$ as a holomorphic function on $\Bers{x_1}$.

Let $\charfunc_{A}$ be the characteristic function of $A$ on $\mathcal{PMF}$. From the property (ii) in the disintegration,
$$
0<\ThursM^{x_0}(A)=
\int_{\mathbb{P}_{x_0}\mathcal{MF}}\left(\int_{\mathcal{PMF}}\charfunc_{A}([F])\,d\lambda_t([F])\right)d\nu^{x_0}(t).
$$
From Proposition \ref{prop:absolutecontiuous}, there are $t\in \mathbb{P}_{x_0}\mathcal{MF}$ and $[F]\in \mathcal{PMF}$ such that $\Pi^{x_0}([F])=t$, $\lambda_t=\Theta_t$ under the identification \eqref{eq:identification_Pi} and
\begin{equation}
\label{eq:radial_F-M-Riesz}
\Theta_t(\{\theta\in\mathbb{S}^1\mid [v(e^{i\theta}q_{F,x_0})]\in A\})>0.
\end{equation}

Consider the Teichm\"uller disk $\mathcal{R}_{[F]}\colon \mathbb{D}\to \teich_{g,m}$ which is defined by $q_{F,x_0}$ with $\mathcal{R}_{[F]}(0)=x_0$. From the assumption, the radial limit of a bounded holomorphic function $f\circ \mathcal{R}_{[F]}$ on $\mathbb{D}$ vanishes at the direction in $A$. From \eqref{eq:radial_F-M-Riesz} and the (original) F. and M. Riesz theorem, we get $f\circ \mathcal{R}_{[F]}\equiv 0$ on $\mathbb{D}$ and hence $f(x_0)=f\circ \mathcal{R}_{[F]}(0)=0$.

\subsection{Proof of Corollary \ref{coro:2}}
Fix a symplectic basis $\{A_i,B_i\}_{i=1}^g$ on $\Sigma_g$ and define the period map $\Pi$ on $\teich_g$. Namely,
for $x=(X,f)\in\teich_g$, let $\psi_i^x$ be the holomorphic $1$-form on $X$ with
$$
\int_{f(A_j)}\psi_i^x=\delta_{ij}\quad (\mbox{Kronecker's delta})
$$
for $1\le i,j\le g$. Let $\displaystyle\pi_{ij}(x)=\int_{f(B_j)}\psi_i^x$ and set $\Pi(x)=\begin{bmatrix}\pi_{ij}(x)\end{bmatrix}$ for $x\in \teich_{g}$.
Then, $\Pi$ is holomorphic on $\teich_g$ and the image of $\Pi$ is contained in the Siegel upper-half space of degree $g$ (cf. \cite{Ahlfors:1960} and \cite{Rauch:1965}). Since the Siegel upper-half space of genus $g$ is biholomorphic to a bounded domain, there is a holomorphic map $\Phi$ defined on the Siegel upper half plane such that all entry of $H:=\Phi\circ \Pi$ is a bounded holomorphic function on $\teich_g$ (e.g. \cite[Theorem 1 in \S3, Chapter 6]{Siegel:1989}). From Theorem \ref{thm:main_radial_limit}, a holomorphic map $H$ admits the radial limits $H^*$ (in our sense). Notice that Shiga \cite[\S5]{Shiga:1984} also discusses the boundary behavior of the period map.

From the definition of $\mathcal{I}_g$, $H\circ [\omega]=H$ on $\teich_g$ for $[\omega]\in \mathcal{I}_g$. 
From \eqref{eq:naturality_MCG}, the radial limit $H^*$ is invariant under the action of $\mathcal{I}_g$ on $\partial^{ue}\Bers{x_0}$. The function $H$ is not constant function since the period map defines local charts at almost every point on $\teich_{g}$ (e.g. \cite{Ahlfors:1960}). Hence, from Theorem \ref{thm:F-M-Riesz}, $H^*$ is also not constant as a (bounded) measurable function on $\partial \Bers{x_0}$. Therefore, a function
$$
\mathcal{PMF}\ni [F]\mapsto
\begin{cases}
H^*\circ \Xi_{x_0}([F]) & ([F]\in \mathcal{PMF}^{ue}) \\
0 & (\mbox{otherwise})
\end{cases}
$$
becomes a non-constant measurable function on $\mathcal{PMF}$ which is invariant under the action of $\mathcal{I}_g$. This implies that the action of $\mathcal{I}_g$ on $\mathcal{PMF}$ is not ergodic.

\section{Conclusion}
In view of Fatou's research \cite{Fatou:1906}, a natural problem next to our result is to present bounded pluriharmonic functions by the Poisson integral. The Poisson integral presentation will characterize the image of the isometry \eqref{eq:isometric_embedding}. Indeed, it is conjectured from the Poisson integral formula in \cite{Miyachi:2023} that the image coincides with the subspace of $L^\infty(\partial \Bers{x_0})$ defined by
$$
\left\{
g\in L^\infty(\partial \Bers{x_0})\mid
\int_{\partial^{mf}\Bers{x_0}}g(\varphi)\overline{\partial}_x
\left\{
\left(
\dfrac{\ext_{x_0}(F_\varphi)}{\ext_{x}(F_\varphi)}
\right)^{3g-3+m}
\right\}
d\pushThursMBers_{x_0}(\varphi)=0\ (x\in \teich_{g,m})
\right\},
$$
where $\overline{\partial}_x$ is the $\overline{\partial}$ derivative in terms of the variable $x\in \teich_{g,m}$, $\partial^{mf} \Bers{x_0}$ is the part of the Bers boundary whose  and $F_\varphi$ is the measured foliation whose singular foliation corresponds to the ending lamination of $\varphi\in \partial^{mf} \Bers{x_0}$ (cf. \S\ref{subsec:BoundarygroupswithoutAPT}).

By taking the pull-back via the map \eqref{eq:curve-complex-mininal-bers}, the image of the isometry \eqref{eq:isometric_embedding} is identified with an invariant closed subspace of $L^\infty(\mathcal{PMF})=L^\infty(\mathcal{PMF},\mu_{Th}^{x_0})$ under the $\mathbb{C}$-linear action of the mapping class group. From the above mentioned conjecture, the space is possibly described as a (closed) subspace of $L^\infty(\mathcal{PMF})$ consisting of $h\in L^\infty(\mathcal{PMF})$ with
$$
\int_{\mathcal{PMF}}h([F])\overline{\partial}_x
\left\{
\left(
\dfrac{\ext_{x_0}(F)}{\ext_{x}(F)}
\right)^{3g-3+m}
\right\}
d\ThursM^{x_0}([F])=0
$$
for all $x\in \teich_{g,m}$.
The closed subspace obtained here reflects the complex structure of the Teichm\"uller space. The $\mathbb{C}$-linear action gives a faithful linear presentation of the mapping class group, except for the finite cases where $(g,m)=(1,1)$, $(0,4)$, $(1,2)$, and $(2,0)$ (e.g. \cite{Papadopoulos:2015}). In the exceptinal cases, the kernel of the action is finite. Thus, the further study of the action is expected to contribute to approach the conjectures mentioned in \S\ref{subsec:background}.

\subsection*{Acknowledgements}
The author thanks Professor Athanase Papadopoulos and Professor Ken'ichi Ohshika for fruitful conversations. The author also thanks Professor Howard Masur for valuable comments. The author thanks anonymous referees for useful comments.



\def\cprime{$'$} \def\cprime{$'$}

\end{document}